%% file: main.tex
\documentclass[runningheads]{llncs}
\usepackage[T1]{fontenc}
\usepackage{graphicx}
\usepackage{xspace}
\usepackage{calrsfs}
\usepackage{color}
\usepackage{amsmath}
\usepackage{amsfonts,amssymb}
\usepackage{url}
\usepackage[overload]{empheq}%for \empheqlbrace
\usepackage{scalerel}
\usepackage{hyperref}
\usepackage{multirow}

\graphicspath{{./figures/}}   %path of figures, always use two {{ and }}

\newcommand{\ie}{\textit{i.e.}}
\newcommand{\eg}{\textit{e.g.}}
\newcommand{\drow}[1]{\multirow{2}{*}{#1}}
\def\cutplanes{\texttt{Cutting-Planes}\xspace}
\def\conb{\texttt{ConicBundle}\xspace}
\def\matlab{\texttt{Matlab}\xspace}
\def\mosek{\texttt{Mosek}\xspace}
\def\pcp{\texttt{Projective Cutting-Planes}\xspace}
\def\projcutplanes{\pcp}

\DeclareMathAlphabet{\pazocal}{OMS}{zplm}{m}{n}
\def\AA{\pazocal{A}}
\def\na{\text{n.a.}}
\def\CC{C}
\def\cols{\mathcal{C}}
\def\dols{\mathcal{D}}
\def\diag{\texttt{diag}}
\def\ub{\texttt{ub}}
\def\lb{\texttt{lb}}
\def\optsol{{\tt opt}}
\def\x{\mathbf y}

\def\zero{\mathbf 0}

\def\itr{\texttt{it}}
\def\itt{\itr}
\def\u{\mathbf u}
\def\v{\mathbf v}
\def\b{\mathbf b}
\def\a{{\mathbf a}}
\def\ca{c_a}
\def\y{\mathbf y}
\def\itr{\texttt{it}}
\def\d{\mathbf d}
\def\p{\mathbf p}
\def\PP{\mathcal{P}}
\def\S{X}
\def\R{\mathbb R}
\def\zeros{{\mathbf{0}}}
\def\ones{{\mathbf{1}}}
\def\lambdamin{\lambda_{\min}}

\def\KncN{[K_{nc}~N]}
\def\midx{X_{c+m}}
\def\midmat{D_{c+m}}
\def\midm{\midmat}
\DeclareMathOperator*{\sprod}{\scalerel*{\cdot}{\bigodot}}

\begin{document}
\title{Semidefinite Programming by {\tt Projective Cutting Planes}}
\author{Daniel Porumbel\inst{1}}
\authorrunning{Daniel Porumbel}

\institute{Conservatoire National des Arts et M\' etiers, Paris, France \\ 
\email{daniel.porumbel@cnam.fr}\\
\url{http://cedric.cnam.fr/~porumbed/}}

\maketitle              
\begin{abstract}
Seeking tighter relaxations of combinatorial optimization problems, 
semidefinite programming is a generalization of linear
programming that offers better bounds 
and is still polynomially solvable. Yet, in practice, a semidefinite program is
still significantly harder to solve than a similar-size Linear Program (LP).
It is well-known that a semidefinite program can be written as an LP with 
infinitely-many cuts that could be solved by repeated separation in 
a \cutplanes scheme; this approach
is likely to end up in failure. We proposed in~\cite{mainpaper} the \pcp
that
%use \cutplanes by 
upgrade the well-known
separation sub-problem to the projection sub-problem: given 
a feasible $\y$ inside a polytope $\PP$ and a direction $\d$, find
the maximum $t^*$ so that $\y+t^*\d\in\PP$. Using this new sub-problem,
one can generate a sequence of both inner and
outer solutions that converge to the optimum over $\PP$. This paper shows that the projection
sub-problem can be solved very efficiently in a semidefinite programming context, enabling
the resulting \pcp to compete very well with state-of-the-art
semidefinite optimization software (refined over decades).
Results suggest it may the fastest method for matrix sizes larger than
$2000\times 2000$.

\end{abstract}
\section{Introduction}

We consider the following semidefinite optimization problem
\begin{subequations}
\label{sdp1}
\begin{align}[left ={(SDP)  \empheqlbrace}]
\max_{\y\in \R^n}~~&\b^\top \y                                                       \label{sdp1a}\\
s.t~~&\AA^\top \y  \preceq C                                            \label{sdp1b}\\
     &\a^\top \y \leq \ca~\forall (\a,\ca)\in \cols                     \label{sdp1c} 
    % &\y\in \R^k,                                                        \label{sdp1d}
    &\end{align}
\end{subequations}
where $\AA^\top \y = \sum_{i=1}^{k} A_iy_i$; $A_1,~A_2,\dots, A_k$ and $C$ are symmetric 
$n\times n$ matrices.
The set $\cols$ in \eqref{sdp1c} contains simple linear constraints that may include $\y\geq 
\zeros$ (\ie, one can enforce $y_i\geq 0$ by taking $a_i=-1$ and $a_j=0\forall j\neq i$ and
 $\ca=0$).
Adding a limited number of other linear constraints in \eqref{sdp1c} would not have a huge impact on 
the techniques we propose, because we will also  
express the semidefinite positive (SDP) constraint \eqref{sdp1b} through a set of linear constraints. 

We rewrite \eqref{sdp1a}-\eqref{sdp1c} in the form below that
is better suited to \cutplanes. Constraint \eqref{sdp1b} is equivalent to
\eqref{sdp2d} when $\dols=\R^n$, \ie, when \eqref{sdp2d} 
incorporates all cuts
$\dols=\R^n$, \eqref{sdp2d} actually
reduces to
$\S\succeq \zeros$ for $\S=  C-\AA^\top \y$.
Recall that
$\S\succeq \zeros \iff \S\sprod \d\d^\top \geq 0~\forall \d\in\R^n$,
where
$\S\sprod \d\d^\top = \d^\top \S \d$. 

\begin{subequations}
\label{sdp2}
\begin{align}[left ={(SDP$-$LP)  \empheqlbrace}]
\max_{\y\in \R^n}~~&\b^\top \y                                           \label{sdp2a}\\
s.t~~& \S=  C-\AA^\top \y                                   \label{sdp2b}\\
     &\a^\top \y \leq \ca~\forall (\a,\ca)\in \cols         \label{sdp2c}\\
    &\S\sprod \d\d^\top \geq 0~\forall \d\in\dols           \label{sdp2d}
    &\end{align}
\end{subequations}

In a \cutplanes scheme,
both sets $\cols$ and $\dols$ could be generated on the fly,
but this work-in-progress paper
only addresses the case of a fixed (often empty) set $\cols$. The main difficulty is to
make $\dols$ grow along the iterations, up to the point of finding the optimum
solution of \eqref{sdp1a}-\eqref{sdp1c} using a reasonably-sized $\dols$
in \eqref{sdp2a}-\eqref{sdp2d}.
We aim at large-scale SDP optimization,
 with a value of $n$ reaching a few thousands and a $k$ reaching a value of hundreds.

Interior Point Methods (IPMs) are a very popular approach
for SDP optimization. These methods can offer everything one
can desire in theory, but may be too slow for $n\geq 2000$, because
they would have to solve huge Newton systems. Each iteration may easily 
end up in requiring computing a positive definite matrix of size $k\times k$
which requires 
$O(k^2n^2+kn^3)$ operations, see, \eg, \cite[p. 66]{helmberg} or \cite[p. 26]{kartik}.
The state of the art solver \mosek implemented one of the fastest IPM for SDP optimization.

Another very successful approach is the
\conb method~\cite{helmbergbase,helmbergsoft} that reformulates
the semidefinite program as an eigenvalue optimization problem, 
which is then solved by a subgradient method. 
This eigenvalue optimization problem arises as follows.
The \conb requires a
constant trace constraint in the dual. Consider the
expression  $\AA^\top \y= \sum_{i=1}^{k} A_iy_i$ 
from \eqref{sdp2b} and suppose we have $A_k=I_n$ and $b_k>0$.
The optimal way to enforce  
$C-\AA^\top \y\succeq \zeros$ and maximize the objective is
to set $y_k=\lambdamin \left(C- \sum_{i=1}^{k-1} A_iy_i\right)$.
Separating this variable from the other $k-1$, one has to maximize this minimum eigenvalue function over the
decision variables $y_1,y_2,\dots, y_{k-1}$.
Such methods have to maintain a cutting model that overestimates
the concave non-smooth minimum eigenvalue function. 
The use of the term $A_k=I_n$ with $b_k>0$ 
reduces to imposing {\tt trace}$(Z)=b_k$ on the
dual matrix $Z$ in the dual of \eqref{sdp2a}-\eqref{sdp2d} --
and one can replace $k$ with a linear combination of $A_1,A_2,\dots A_k$.

The \pcp method proposed in this paper was deliberately designed to be as lightweight
as possible, even more so than the bundle methods. 
Paradoxically, this is both a strength and a weakness. The weakness is  
that most proposed ideas are rather ad-hoc and they do not emerge in a structured
manner from an established theory that has the size or strength
to attract many people; this explains why the paper has 
few bibliographical references. The advantage is that the new method
is not placed on an existing research thread 
acknowledged for a long time in SDP optimization but (only) followed by highly-specialised 
SDP professionals. The number
of people who can understand the new approach is rather large, extending (a bit)
beyond the SDP field.

\section{The general \pcp}

Let $\PP$ be the feasible area of the semi-infinite LP \eqref{sdp2a}-\eqref{sdp2d}.
A \cutplanes algorithm constructs at 
each iteration $\itr$ an outer
approximation $\PP_\itr$ of $\PP$, \ie, a polytope
$\PP_\itr$ defined only by
 a subset of the constraints of  $\PP$, so that
$\PP_\itr\supseteq \PP$. 
This generates a sequence of upper
bounds $\b^\top \optsol(\PP_\itr)$ that decrease along the
iterations $\itt$, converging to 
the optimal solution value 
$\b^\top \optsol(\PP)$; these bounds are
associated with a series of 
outer non-feasible
solutions $\optsol(\PP_\itr)$.

Any outer solution $\optsol(\PP_\itr)$ of \eqref{sdp2a}--\eqref{sdp2d}
obtained this way can be turned into a feasible solution $Z$ of the dual
of the main SDP program
\eqref{sdp1a}-\eqref{sdp1c}. It is enough to restrict $\dols$ in 
\eqref{sdp2d} to the current set $\overline{\dols}$ of 
active constraints. In fact, each $\d\in\dols$ 
generates a constraint \eqref{sdp2d} has to be understood in conjunction with
\eqref{sdp2b} and implemented as $
\left(A_1\sprod \d\d^\top\right) y_1+
\left(A_2\sprod \d\d^\top\right) y_2+
\dots
\left(A_k\sprod \d\d^\top\right) y_k
\leq
C\sprod \d\d^\top 
$.
By LP duality, the objective value of the outer
LP solution is the same as that of the dual LP solution 
$Z=\sum_{\d\in\overline{\dols}}\gamma_\d \d\d^\top\succeq\zeros$, where $\gamma_\d$
is the optimal dual value of constraint \eqref{sdp2d} for
$\d$; this $Z$ is also the dual SDP solution mentioned above.
See also \cite[Theorem 9]{kartik} for more details on how a
\cutplanes solving \eqref{sdp2a}-\eqref{sdp2d} can generate
feasible dual SDP solutions along the iterations.
But there is no 
built-in functionality in \cutplanes to generate 
inner feasible solutions.

\begin{figure}[ht]
\hspace{-1.0em}
\scalebox{0.82}{%
\ifx\pdfoutput\undefined
  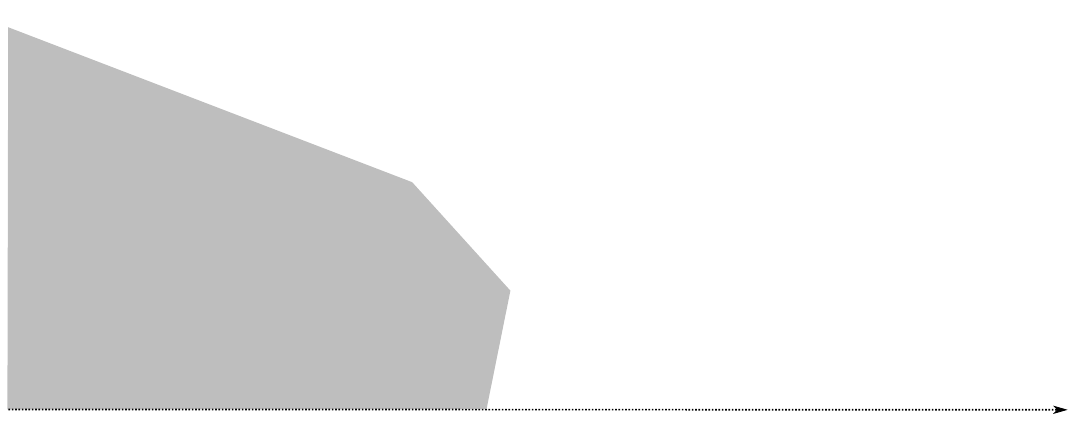
\else
  \ifnum\pdfoutput=1 %this might be enough
    \input{./figures/fig1.pdf_tex}
  \else
    \input{./figures/fig1.pdf_tex}
  \fi
\fi
}
\caption{\small \label{fig1}Example of \projcutplanes execution (3 iterations).  At the
first iteration, one may projects $\x_1=\zero$ towards the optimal solution of an initial 
default polytope $\PP_0$ that may contain only (very loose) bounds on $\y$. The projection sub-problem returns $t_1^*$ and the 
first-hit constraint represented by the black exterior solid line.  At
iteration $\itt=2$, if we use $\alpha=0.5$, the midpoint $\x_2$ between $\x_1$ and 
$\x_1+t^*\d_1$ is projected towards 
the optimal outer solution $\optsol(\PP_1)$ ---
at iteration 1, 
the outer approximation
$\PP_1\supset\PP$ 
contains the largest triangle. This generates a second facet 
(blue solid line) that is added to the facets
of $\PP_1$ to construct $\PP_2$.  The third projection
in red
takes the midpoint 
$\x_3$
between the blue square and the blue circle (last pierce point) and projects it towards $\optsol(\PP_2)$.
}
\end{figure}

To generate both inner and outer solutions (with regards to $\PP$),
\pcp uses an iterative operation
of projecting an interior
point inside $\PP$, as illustrated in Fig.~\ref{fig1}. At each iteration $\itt$, an inner solution $\x_\itt\in\PP$
is projected towards the direction $\d_\itt$
of the current optimal outer solution $\optsol(\PP_{\itt-1})$, \ie, 
we take
$\d_\itt= \optsol(\PP_{\itt-1})-\x_\itt$.
The projection 
sub-problem asks to determine
$t^*_\itt=\max\left\{t:\x_\itt+t\d_\itt\in\PP\right\}$.
This requires
finding
the pierce (hit) point 
$\x_\itt+t^*_\itt\d_\itt$ and a (first-hit) constraint of $\PP$,
which is added to the constraints of $\PP_{\itt-1}$
to construct $\PP_\itt$. At next iteration $\itt+1$,
\projcutplanes takes a new interior 
point $\x_{\itt+1}$ on the segment joining $\x_\itt$ and $\x_\itt+t^*_\itt\d_\itt$ and 
projects it towards direction
$\d_{\itt+1}=\optsol(\PP_\itt)-\x_{\itt+1}$. 
Regarding the theoretical convergence proof, 
see Remark~\ref{remConv} (p.~\pageref{remConv}) in appendix.

A central question in practice is to choose the projection base:
given
$\x_\itt$ and $\x_\itt+t^*_\itt\d_\itt$, how should 
one choose the point $\x_{\itt+1}=\x_\itt+\alpha\cdot t^*_\itt\d_\itt$?
Using $\alpha=1$ would make \pcp very aggressive. But 
previous work on combinatorial optimization LP relaxations (see, \eg, Sections 2.2.2 or 3.2.1
of~\cite{followup}) show that such a choice may only
produce better
feasible solutions in the beginning, but need more iterations in the long
run. Such variant may be, however, useful if we do not need a very tight 
gap ${\tt ub}_\itt-{\tt lb}_\itt=\b^\top \optsol(\PP_\itr) - \b^\top \y_\itt$; for instance, if we solve a relaxation of a combinatorial
optimization problem, we can stop when 
$\lfloor {\tt ub}_\itt \rfloor=
 \lfloor {\tt lb}_\itt \rfloor$.
Based on previous work, we decided to use
a rather conservative
step length $\alpha=0.3$; better choices may exist. 
In this work-in-progress paper, we start at the very first iteration with
$\x_1=\zeros_k$, but we have already studied other options that will be submitted
for publication in a longer paper.

To determine
$t_\itt^*=\max\left\{t:\x_\itt+t\d_\itt\in\PP\right\}$,
one also has 
to find a {\it first-hit constraint} satisfied with equality by
$\x_\itt+t_\itt^*\d_\itt$. This projection sub-problem implicitly solves the separation sub-problem
for all points  $\x_\itt+t\d_\itt$ with $t\in\R_+$, because the above first-hit
constraint separates  all solutions $\x_\itt+t\d_\itt$ with
$t>t_\itt^*$ and proves $\x_\itt+t\d_\itt\in\PP~\forall t\in[0,t_\itt^*]$.

By generalizing the separation sub-problem,
the projection sub-problem
may seem 
computationally far more expensive, 
but we will see this is not necessarily the case.
Section \ref{sec3} 
presents a few SDP techniques
that bring us very close to designing a projection algorithm as fast as
the separation one. Numerical experiments (Section~\ref{sec4}) confirm that,
in general, 
the projection sub-problem is not the most important computational
bottleneck of the overall method.

The new method is reminiscent of an Interior Point Method
(IPM) by the way it generates a sequence of interior points that
converge to
the optimal solution. An IPM 
moves from solution to solution by advancing along a Newton direction at each
iteration, in an attempt
to solve first order optimality conditions~\cite{gondzio25years}. 
Advancing along a Newton 
direction is not really equivalent to performing
a projection, because a projection
advances up to the pierce point,
while a Newton step in an IPM does not even execute all iterations to fully 
solve the first order conditions (for the current barrier term). 
An IPM tries to 
generate
well-centered dual solutions that stay in 
the proximity of a central path; this 
is reminiscent of the trajectory of feasible
solutions $\x_1,~\x_2,~\x_3,\dots$ constructed 
by \pcp.

\par\sloppy 
The remaining (few) SDP customizations needed to adapt \pcp to
an SDP context are listed in
Remark~\ref{remNutsAndBolts} 
(p.~\pageref{remNutsAndBolts} in appendix).
\par

Before presenting the SDP projection (Section~\ref{sec3}), recall
the semi-infinite LP \eqref{sdp2a}-\eqref{sdp2d} may contain two sets
of constraints that may be generated on the fly: \eqref{sdp2c}
end \eqref{sdp2d}. When necessary, one may have to solve two
projection sub-problems, a non-SDP one with regards to \eqref{sdp2c}
and an SDP one with regards to \eqref{sdp2d}. The space limitation 
does not enable us to advance more on this idea, but, to our knowledge, such
questions are out of reach for other SDP algorithms.

%Remark~\ref{remNutsAndBolts} (p.~\pageref{remNutsAndBolts}) in appendix presents in more detail
%all SDP customizations needed to run \pcp in the above LP framework.

\section{The SDP projection algorithm\label{sec3}}

The \pcp was initially designed and tested independently of any SDP 
concept, for the purpose of solving LP relaxations in combinatorial
optimization.
This is the first time we solve the projection sub-problem over the SDP cone:
what is the maximum 
$t^*$ so that $X+t^*D\succeq \zeros$? In our context, $\S$
is the SDP matrix $\S=C-\AA^\top \y_{in}\succeq \zeros$ associated to 
the current inner point $\y_{in}$ of \eqref{sdp2a}-\eqref{sdp2d}
and $D=-\AA^\top (\y_{out}-\y_{in})$, where $\y_{out}$ is the current outer point. We also have to determine a first-hit
vector $\v\in\R^n$ so that $\left( X+t^*D\right)\sprod \v\v^\top =0$. In practice, we may
easily encounter numerical problems and this equality will always be 
satisfied within a certain tolerance. 
On the other hand, the value $D\sprod\v\v^\top $ should be really significantly
lower than $0$. When this is the case,  advancing any
$\epsilon>0$ beyond $t^*$ leads to $\left( X+(t^*+\epsilon)D\right)\sprod\v\v^\top <0$. 

\begin{property}\label{prop1}
We will see that the projection $X\to D$ can be calculated more rapidly if
$D$ belongs to the image of $X$. This means that each column 
(and row) of $D$ can be written as a linear combination of the columns
(or rows, resp.) of $X$. We can equivalently say that the null space of $X$
is included in the null space of $D$; thus, $X\d=0\implies D\d=0~\forall \d\in\R^n$. 
We will show below in cases A) and B) how it is easier to project 
when this property holds; if possible, \pcp should thus adapt its own evolution to 
seek this property.
\end{property}

Two matrices $X$ and $X'$ are congruent if there is some non-singular $M$ 
such that $X'=MXM^\top$. 
It is well known (see, for
example, \cite[Prop 1.2.3.]{mysdp}) that two
congruent matrices have the same SDP status: $X\succeq \zeros
\iff X'\succeq \zeros$.

\begin{property} (congruent expansion)\label{defexpansion}
We say that
$X'\in\R^{n'\times n'}$ with $n'>n$
is a congruent expansion of
$X\in\R^{n\times n}$ if and only if we can write $X'=MXM^\top$,
for some $M\in \R^{n'\times n}$ of {\it full rank $n$}. 
$X$ has the same SDP status as $X'$.

\end{property}
\begin{proof}We show both implications below.
\begin{enumerate}
\item $X\succeq\zeros \implies X'\succeq \zeros$. 
Assume the contrary for the sake of contradiction: $\exists \v'\in\R^{n'}$ such 
that $\v'^\top X'\v'<0$. This implies 
$\v'^\top MXM^\top\v'<0$, equivalent to $X\nsucceq\zeros$, contradiction.
\item $X'\succeq\zeros \implies X\succeq \zeros$
Assume the contrary: $\exists \v\in\R^{n}$ such 
$\v^\top X\v<0$.
We can surely write $\v^\top = \v'^\top M$ for some $\v'\in\R^{n'}$ because
$M$ has full rank. This means $\v'^\top M X M^\top \v'<0$, equivalent to $\v'^\top X'\v'<0$,
contradiction. \qed
\end{enumerate}
\end{proof}

Using these concepts we are ready to address the projection algorithm. We will distinguish
four cases noted A), B), C) and D). Given $X\succeq\zeros$, we have to find:
\begin{equation}\label{eqmaxt} \max\left\{t:X+tD\succeq \zeros\right\}.\end{equation}

\noindent \line(1,0){33}\\
\noindent A) This case is characterized  by $X\succ\zeros$, \ie, $X$ is non-singular;
Prop~\ref{prop1} surely holds because the image of a non-singular
$X$ is $\R^n$. We apply the Cholesky decomposition to determine
the unique non-singular $K$ such that $X=KK^\top$.
We then solve
$D=KD'K^\top$ in variables $D'$ by back substitution; this may require $O(n^3)$ in theory, 
but \matlab is able to compute it much more rapidly in practice
because $K$ is triangular.
Let us re-write \eqref{eqmaxt} as:

\begin{equation}\label{eqmaxt2}
\max\left\{t:KI_nK^\top+tKD'K^\top\succeq \zeros\right\}.
\end{equation}
This is equivalent (by congruence according to Prop~\ref{defexpansion}) to
\begin{equation}\label{eqmaxt3}
\max\left\{t:I_n+tD'\succeq \zeros\right\}.
\end{equation}
The sought step length is $t^*=-\frac 
1{\lambdamin(D')}$, or $t^*=\infty$ if $\lambdamin(D')\geq 0$.

We still have to find a first-hit cut $\v\in\R^n$; in fact, technically,
the first-hit cut will be
$
\left(A_1\sprod \v\v^\top\right) y_1+
\left(A_2\sprod \v\v^\top\right) y_2+
\dots
+
\left(A_k\sprod \v\v^\top\right) y_k
\leq
C\sprod \v\v^\top 
$.

If $\v$ is an eigenvector of $K (I_n+t^*D') K^\top$ with an eigenvalue of 0, this means
$\v^\top K (I_n+t^*D') K^\top \v=0$. Thus, $\u=K^\top \v$ is eigenvector
of $I_n+t^*D'$ with an eigenvalue of 0. This latter eigenvector $\u$ can be computed when determining
$\lambdamin(D')<0$ above, because if the eigenvalue of $\u$ with regards to $D'$ is
$\lambdamin(D')$ its eigenvalue with regards to $I_n+t^*D'$ is 0 (since recall
$t^*=-\frac 
1{\lambdamin(D')}$).
The sought $\v$ solves $K^\top \v=\u$ and it can rapidly be computed
by back-substitution. We have $\u^\top D'\u<0 \implies
\v^\top K D' K^\top \v<0 \implies \v^\top D \v<0$. We thus
have 
$\v^\top (X+t^*D) \v =0$
and
$\v^\top (X+(t^*+\epsilon)D) \v <0$ for any $\epsilon>0$. This proves
$\v$ is a first-hit cut.

%The approach will be generalized for other matrices besides $X+t^*D$.

\noindent \line(1,0){33}\\
\noindent B) In this case Prop~\ref{prop1} is still satisfied, but $X$ has rank $c<n$.
This means $X$ contains
$c$ independent rows (and columns by symmetry), referred to 
as {\it core} rows (or columns); the other dependent rows (or 
columns) are 
{\it non-core} positions.
%We can 
%always rearrange the rows and columns of $X$ so that the top-left matrix of 
%size $c\times c$ is full-rank and strictly SDP.
%And we can apply the same operation on $D$. 
Using the LDL decomposition of $X$, we will factorize
$X=K_{nc}K_{nc}^\top$, where $K_{nc}\in\R^{n\times c}$.
The image of $K_{nc}$ is equal to the image of $X$. 
Since 
Prop~\ref{prop1} is satisfied, we will see we can still solve
$D=K_{nc}D'K_{nc}^\top$ in variables $D'$. A first intuition is to notice that we 
can project $X\to D$ only over the core rows and columns, 
because the non-core positions are dependent on the core 
ones.

But the most difficult task is to determine these core positions.
We first apply the LDL decomposition and write $X=L\diag(\p)L^\top$ with
$\p\geq \zeros_n$.
The contribution of each $p_i$
in $L\diag(\p)L^\top$ is actually $p_i L_iL_i^\top$, where
$L_i$ is column $i$ of $L~(\forall i\in[1..n])$. If all $n\times n$ elements of
$p_{i}L_iL_i^\top$ are below some precision parameter,
we consider $i$ is a non-core position; otherwise, it is a core position.
By reducing all non-core positions
$p_{i}$ to zero, we can say that all $n-c$ non core
columns of $L$ vanish in the decomposition
$X=L\diag(\p)L^\top$. After 
removing these vanished $n-c$ columns from $L$ and the corresponding 
zeros from $\p$, 
we can write
$X=L\diag(\p)L^\top=L\diag(\p)^\frac12\diag(\p)^\frac12L^\top = K_{nc}K_{nc}^\top$ with $K_{nc}\in\R^{n\times c}$.

We next solve
$D=K_{nc}D'K_{nc}^\top$ in variables $D'$. For this, we first reduce
this system to work on $c\times c$ matrices, \ie, we transform
it into
$D_{cc}=K_{cc}D'K_{cc}$ where 
$K_{cc}$ is $K_{nc}$ restricted to the $c$ core rows
and
$D_{cc}$ is $D$ restricted to the $c\times c$ core rows and 
columns. To solve this square system, 
we  
apply back-substitution twice and this is very fast because
$K_{cc}$ is lower triangular. If the resulting solution $D'$ also
satisfies $D=K_{nc}D'K_{nc}^\top$,
then we are surely in case B).
We obtained a reduced-size version of
\eqref{eqmaxt3} working in the space of $c\times c$ matrices:

\begin{equation}\label{eqmaxtb}
\max\left\{t:I_c+tD'\succeq \zeros\right\}.
\end{equation}
And the maximum value of $t$ is here: $t^*=-\frac 
1{\lambdamin(D')}$, or $t^*=\infty$ if $\lambdamin(D')\geq 0$.

We finally determine a first-hit vector 
$\v_c\in \R^c$ 
over the core rows and columns exactly like in (the last paragraph describing) case A).
To lift $\v_c$ to a 
hit vector $\v\in\R^n$, we construct $\v$ by inheriting
the core positions from $\v_c$ and filling the non-core
positions with zeros.

\noindent \line(1,0){33}\\
\noindent C) We still use the decomposition 
$X=K_{nc}K_{nc}^\top$ computed above but we suppose
that the system $D=K_{nc}D'K_{nc}^\top$ has no solution in variables $D'$.
This also means Prop~\ref{prop1} is not satisfied:
$D$ does not belong to the image of $K_{nc}$ or $\S$.

We will express all columns of $D$ as a linear
combination of: (i) the columns of $K_{nc}$ and (ii) a set 
of $m$ columns of
$D$ named active (independent) columns. We first apply the QR decomposition on matrix 
$
\left[K_{nc}~D \right]\in \R^{n\times (c+n)}$ and write
$\left[K_{nc}~D \right] = QR$, where $Q\in\R^{n\times (c+n)}$
and $R\in\R^{(c+n)\times (c+n)}$
is upper triangular.
In fact, the standard QR factorization returns 
a matrix $Q\in\R^{n\times n}$
and a matrix $R\in\R^{n\times (c+n)}$, but we artificially extend
$Q$ with $c$ null columns and $R$ with $c$ null rows to simplify 
notations.
Let us focus on the first $c$ columns of 
$R$. Since $K_{nc}$ is full rank, the matrix $R$ restricted 
to the first $c$ columns will be full rank; since it is upper 
triangular, this means $R_{jj}\neq 0$ for all $j\in[1..c]$.

Now focus on row $c+i$ of $R$ for each $i\in[1..n]$. If all elements of
this row are zero, column $c+i$ of $Q$ it has no contribution in the product $QR$; this also means 
that column $c+i$ of  $\left[K_{nc}~D \right]$ can be 
expressed as a combination of the first $c+i-1$ columns of $Q$. We call
this column of $Q$ non-active, being dependent on
the columns of $K_{nc}$ and on the active
columns found while scanning the columns $[c+1,c+i-1]$ of $Q$.
%The QR decomposition algorithm might have very well ignored 
%constructing column $c+i$
%of $Q$ in such case, because it has no impact in the product $QR$.

Let $N$ denote the matrix $Q$ restricted to its $m$ active 
columns detected above. The size of $N$ provides a new way to detect  
case B): if $N$ were empty with $m=0$, we would have been in case B). When $N$ is non-empty, we can decompose
$X$ and $D$ as follows:
\begin{align}\label{decompx}
X&=
\underbrace{
\begin{bmatrix}
        K_{nc}     &    N     \\
\end{bmatrix}
}_{c~+m}
\begin{bmatrix}
        I_c  ~~   &     0  ~   \\
        0    ~ &     0~
\end{bmatrix}
\begin{bmatrix}
        K_{nc}^\top\\
        N^\top
\end{bmatrix}\\
\label{decompd}
D&=
\underbrace{
\begin{bmatrix}
        K_{nc}     &    N     \\
\end{bmatrix}
}_{c~+~m}
\underbrace{
\begin{bmatrix}
        F     &    G^\top     \\
        G     &     E
\end{bmatrix}
}_{\midmat}
\begin{bmatrix}
        K_{nc}^\top\\
        N^\top
\end{bmatrix}\\
&=K_{nc}FK_{nc}^\top+NEN^\top+K_{nc}G^\top N^\top+NGK_{nc}^\top
\label{decompd2}
\end{align}

The hardest computational task is computing
$\midmat$. A straightforward approach may be quite
slow. We prefer to exploit again the information determined
by the QR decomposition. 
We will modify both sides of the factorization $\left[K_{nc}~D \right] = QR$ to make it similar to \eqref{decompd}.
To transform $Q$ into $\KncN$, we write $Q=[Q_{nc}~Q_{n,c+1..n}]$, \ie, we split its first 
$c$ columns 
$Q_{nc}$ from the last $n$ columns $Q_{n,c+1..n}$.
We can write $K_{nc}=Q_{nc}R_{cc}$, where $R_{cc}$ is the $c\times c$ top-left part of
$R$. Since this system is full rank, we obtain
$Q_{nc}= K_{nc}R_{cc}^{-1}$. We can thus write
$Q=
[K_{nc}~Q_{n,c+1..c+n}]
\left[
\begin{smallmatrix}
        R_{cc}^{-1} & \zeros\\
   \zeros           & I_n       
\end{smallmatrix}
\right]$.
Replacing this $Q$ in $\left[K_{nc}~D \right] = QR$, we obtain
$\left[K_{nc}~D \right] = 
[K_{nc}~Q_{n,c+1..n}]
\left[
\begin{smallmatrix}
        R_{cc}^{-1} & \zeros\\
   \zeros   & I_n       
\end{smallmatrix}
\right] R.$ We now compute the last $n$ columns of the right multiplication and
denote the result by $\CC\in\R^{(n+c)\times n}$. If we 
also restrict $\left[K_{nc}~D \right]$ to its
last $n$ columns (\ie, to $D$), the above QR factorization becomes:

\begin{equation}\label{eq9} D = [K_{nc}~Q_{n,c+1..c+n}] \CC \end{equation}

The matrix $N$ is simply $Q_{n,c+1..c+n}$ restricted to the $m$ active 
columns identified above -- recall a non-active column $c+i$ of $Q$ has
no contribution in the $QR$ product since row $i$ of $R$ is null. The left factor
$[K_{nc}~Q_{n,c+1..c+n}]$  in \eqref{eq9} is thus reduced to $\KncN$ by removing 
the non-active columns of $Q$. The associated $\CC$ factor in 
\eqref{eq9} is also reduced to some $\overline{\CC}$ by removing its null rows that
come from the null rows of $R$.
We can re-write \eqref{eq9} as 
$D = \KncN \cdot \overline{\CC}$. We finally
determine $\midmat$ from \eqref{decompd} by solving
$\midmat 
\left[
\begin{smallmatrix}
        K_{nc}^\top\\
        N^\top
\end{smallmatrix}\right] = \overline{\CC}$. We recall
$\midmat$ has the form:

$$
\midmat = 
\begin{bmatrix}
        F     &    G^\top     \\
        G     &     E
\end{bmatrix}
$$

The case C) under discussion here is characterized by the fact that $G=\zeros$. 
Applying~\eqref{decompd2}, this means
$D$ has the form $D=K_{nc}FK_{nc}^\top+NEN^\top$, where $N$ is orthogonal to 
$K_{nc}$ by (the QR decomposition) construction.\footnote{As a side remark, we can show that $XD$ belongs
to the column image of $X$ in this case C). Since $XD=K_{nc}K_{nc}^\top (K_{nc}FK_{nc}^\top+NEN^\top)$
and $K_{nc}^\top N=0$, we have $X\d=0\implies K_{nc}K_{nc}^\top \d=0 \implies 
K_{nc}^\top \d=0 \implies
XD\d=0,~\forall \d\in\R^n$.}
Using the congruence expansion Property~\ref{defexpansion} on 
\eqref{decompx}-\eqref{decompd}, the SDP status of $X+t\cdot D$ is the
same as that of 
\begin{equation}\label{eqinner}
\begin{bmatrix}
        I_c  ~~   &     0  ~   \\
        0    ~ &     0~
\end{bmatrix}
+t\cdot
\begin{bmatrix}
        F     &    \zeros     \\
        \zeros     &     E
\end{bmatrix}.
\end{equation}
Any hit-vector $\v'\in\R^{c+m}$ for the projection problem 
$\left[
\begin{smallmatrix}
        I_c  ~~   &     0  ~   \\
        0    ~ &     0~
\end{smallmatrix}\right]
\to
\left[\begin{smallmatrix}
        F     &    \zeros     \\
        \zeros     &     E
\end{smallmatrix}\right]
$ can be lifted to a hit-vector $\v\in \R^n$ for the original
projection $X\to D$ by finding a solution $\v$ of the underdetermined
system
$\v'= 
\left[
\begin{smallmatrix}
        K_{nc}^\top\\
        N^\top
\end{smallmatrix}
\right]\v$. We can 
hereafter only focus on projecting 
$\left[
\begin{smallmatrix}
        I_c  ~~   &     0  ~   \\
        0    ~ &     0~
\end{smallmatrix}\right]
\to
\left[\begin{smallmatrix}
        F     &    \zeros     \\
        \zeros     &     E
\end{smallmatrix}\right]$. Case C) is split in two cases:
\begin{description}
\item[C.1)]  if $E\succeq\zeros$, the above projection reduces to 
determining $\max\left\{t:I_c+tF\succeq \zeros\right\} $ and this is solved
using \eqref{eqmaxt} as in case A).
\item[C.2)] 
if $E\nsucceq\zeros$, the sought $t^*$ is 0. It is straightforward to see in 
\eqref{eqinner} that any $t>0$ would generate in this case a non SDP matrix.
\end{description} 

\noindent \line(1,0){33}\\
\noindent D) If all above cases fail, we still 
apply the logic of \eqref{eqinner} and
solve the projection by
finding the maximum $t$ such that:

\begin{equation}\label{eqinner2}
\underbrace{
\begin{bmatrix}
        I_c  ~~   &     0  ~   \\
        0    ~ &     0~
\end{bmatrix}
}_
{\midx\in\R^{(c+m)\times(c+m)}}
+t\cdot
\underbrace{
\begin{bmatrix}
        F     &    G     \\
        G     &     E
\end{bmatrix}}_{\midmat\in\R^{(c+m)\times(c+m)}} \succeq \zeros.
\end{equation}

We first present a tricky case. If
there is some $i\in[c+1..m]$ and some $j\in[1..c]$ such that
the diagonal element $(i,i)$ of $\midmat$ is zero while its non-diagonal element $(i,j)$ is 
non-zero, then
any $t>0$ leads to $\midx+t\midmat\nsucceq \zeros$. We return $t^*=0$, but there is no hit 
vector $\v$ such that $(\midx+t\midmat)\sprod \v\v^\top <0\forall t>0$.
Because if we reduce the whole projection to rows and columns $i$ and $j$, 
there
is no vector $\v\in\R^2$ such that $\left[\begin{smallmatrix} 1 & t\\ 
t & 0 \end{smallmatrix} \right] \sprod \v\v^\top <0$
for any $t>0$ no matter how small.

If all possibilities discussed up to here fail, we solve the projection in two steps:
(1) find a small $t_1$ such that $\midx+t_1\midmat\succeq \zeros$ and (2) solve the projection
$(\midx+t_1\midmat)\to \midmat$. In this second step, $\midmat$ belongs to the image of 
$\midx+t_1\midmat$ 
(Prop~\ref{prop1} satisfied) and we will
use case A) or B). However, finding $t_1$ may require 
a limited number of (costly) repeated separations.
This case is virtually never needed in the experiments
presented in this paper and we explore it further in appendix
(Remark~\ref{remCaseD}, p.~\pageref{remCaseD}).

\section{Numerical results\label{sec4}}
There is unfortunately no well-established benchmark for testing SDP algorithms and no
universally-accepted methodology to measure their performance.
%There is actually no unified benchmark (much less than in fields like vehicle routing or timetabling) 
Most testing has been carried out in rather disparate contexts.
Since we consider the most general SDP programs (no sparsity 
and no particular combinatorial
structure behind the involved matrices), we simply generated 
the instances as follows.
First, we
constructed $\frac k2$ eigenvectors meant to become $0-$eigenvalue eigenvectors (with an eigenvalue of $0$)
for the matrices $A_1,A_2,\dots, A_k$ and $C$; each such eigenvector is inserted in each of these
matrices with a probability of $0.8$.
Once $n_0$ such $0-$eigenvalue eigenvectors are fixed for a given matrix, we construct at random $n-n_0$ orthogonal 
eigenvectors (that together constitute a basis of $\R^n$). In a first 
instance set, we generate the eigenvalues of these $n-n_0$ eigenvectors randomly between $9$ and $10$ for
$A_1,A_2,\dots, A_k $ and between $30$ and $50$ for $C$. In a second instance set, these
eigenvalues have larger variations (indicated by Column 3 of Table~\ref{tab2}). We set $\b=\ones$
when not stated otherwise.

\begin{figure}[t!]
\includegraphics[width=\textwidth]{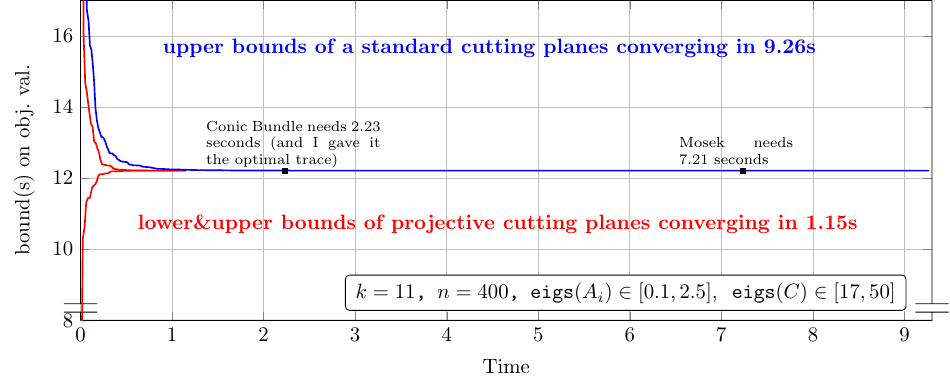}
\vspace{-2em}
\caption{A sample run comparing the main software considered in this paper. 
\label{fig2}}
\end{figure}

Figure~\ref{fig2} illustrates a comparison
between the new method, the standard \cutplanes, the \conb and the \mosek 
solver.
This figure confirms the standard \cutplanes is too slow.
{\tt Mosek} is not very fast for such a low $k$ and large $n$. The \conb needs a bit more than 2 seconds, around twice as
much as \pcp. In this paper we stop \pcp when the {\tt ub-lb} gap is below 0.00001, but
notice that after 0.33 seconds this gap was already hardly noticeable on this figure.
A rather loose gap may be satisfactory when we solve a relaxation of
    combinatorial optimization problem that has an integer optimum.
Since the \conb reformulates \eqref{sdp2a}-\eqref{sdp2d} as an
eigenvalue optimization problem, it needs as input the trace of the optimal
dual solution of \eqref{sdp2a}-\eqref{sdp2d}; we offered it this artificial advantage
by inserting matrix
$A_{k+1}=I_n$; see full details in appendix (Remark~\ref{rem3}, p.~\pageref{rem3}).

\renewcommand{\arraystretch}{1.2}
\begin{table}[ht]
\begin{center}
\begin{tabular}{c|c c|c|c|c|c|c}
\multicolumn{2}{c}{}   &   \multicolumn{6}{c}{n}\\
   \cline{3-8}
\multicolumn{2}{c}{} &  400   & 600     &   800       &   1000  & 1200   & $\stackrel{2000}{{\tiny \text{(gap $ub$-$lb$)}}}$ \\

         %&           &        &         &             &         &         &\\
   %\cline{3-7}
 \multirow{6}{*}{k}
& 20 &$\frac{2.8 (0.7/0.6/0.3)}{9.5/10.6}$  & $\frac{7.1 (2.4/1.6/0.4)}{18.2/34}$ & $\frac{11 (4/2.4/0.3)}{57.5/80}$ & $\frac{26 (10/7/0.5)}{125/134 }$     &$\frac{48   (18/13/0.5)}{220/265 }$   & $\frac{\text{gap closed}}{\na/0.044      }$  \\
& 30 &$\frac{7.8 (2.6/1.6/0.9)}{26.8/15.3}$ & $\frac{22 (9/3/0.2)}{119/43}$       & $\frac{46  (20/8/2.4)}{134/97}$  & $\frac{55 (25/10/1.1)}{155/211}$     &$\frac{125  (55/23/2)}{156/405}$      & $\frac{\text{gap closed}}{\na/0.16       }$  \\
& 40 &$\frac{18  (6 /2/4 )}{108/20}$        & $\frac{26 (10/3/2)} {127/57}$       & $\frac{65(30/8/4 )}{307/134}$    & $\frac{169 (83/21/6)}{391/267}$      &$\frac{253  (123/34 /5 )}{664/479}$   & $\frac{\text{gap closed}}{\na/0.27      }$  \\
& 50 &$\frac{69  (20 /5/25)}{130/25}$       & $\frac{65  (30 /8/3 )}{160/80}$     & $\frac{189(88/17/23)}{410/169}$  & $\frac{297 (151/31 /15)}{863/391}$   &$\frac{462  (237/49 /15)}{770/712}$   & $\frac{0.00023            }{\na/0.94       }$  \\
& 60 &$\frac{76  (22 /4/28)}{128/28}$       & $\frac{178 (68/11/47)}{136/85}$     & $\frac{239(112/19/28)}{658/403}$ & $\frac{520 (270/53/37)}{1081/736}$   &$\frac{614  (318/60 /23)}{1244/975}$  & $\frac{0.0021             }{\na/5.05       }$  \\
& 70 &$\frac{106 (34 /5/38)}{150/35}$       & $\frac{212 (70/11/72)}{474/102}$    & $\frac{628(248/33/187)}{411/735}$& $\frac{570 (304/45/45)}{2283/842}$   &$\frac{1297 (674/103/87)}{2669/1152}$ & $\frac{0.017              }{\na/25.5       }$  \\
\end{tabular}
\end{center}
\caption{ \pcp compared to \conb~\cite{helmbergsoft} and {\tt Mosek}. For each instance, we 
provide the total wall running time (seconds) under the form $\frac{p (a/b/c)}{e/f}$, where
$p$ is the total of wall time of \pcp, $a$ is the time of computing $X$ and $D$, $b$ 
counts the projection time, $c$ is the time of the LP solver for the outer 
approximation of \eqref{sdp2a}-\eqref{sdp2d};
$e$ is the total \conb time 
 and $f$ it the total \mosek time.
 The last column concerns a huge instance size and is different: it provides for each 
 algorithm the gap $\ub-\lb$ reported
after 1000 seconds. In fact, we only report this for our method (the numerator) and \mosek;
``\na'' means
 not available for the \conb, because it does not compute such intermediate $\ub-\lb$ values.
\label{tab1}\vspace{-2em}}
\end{table}

Table~\ref{tab1} reports the wall running times of \pcp, \conb and \mosek on the first instance set.
These are the most time-consuming operations observed on our standard laptop 
(described by Remark~\ref{remComputer}).
\begin{enumerate}
\item[(a)] Determine $X$ and $D$ at each iteration. This operation has complexity 
$O(kn^2)$ while many calculations of the projection algorithm have a complexity 
of $O(n^3)$. Yet these latter calculations use very strongly-optimized \matlab 
routines, while computing $X$ or $D$ can not benefit from such routines, since
this is not a very classical matrix operation. We are almost certain we will
improve this situation in future versions of the software.
\item[(b)] Solve the projection sub-problem $X\to D$. Table~\ref{tab1} show that this
may often represent less than 10\% of the total running time. To our surprise,
the operation from Point (a) is often more computationally expensive.
%This comes from the fact that it works with $k$ matrices of size $n\times n$, while
%the projection mainly works with two matrices of this size ($X$ and $D$).
\item[(c)] Solve the LP corresponding to the outer approximation of the feasible SDP area
\eqref{sdp2a}-\eqref{sdp2d}. This step is relatively insignificant for $k<50$, but
it becomes expensive as $k$ is increased towards 100. The speed
of \pcp for a (much) larger $k$ is dependent on the LP solver (\texttt{cplex});
any future progress in linear programming may bring positive consequences.
\end{enumerate}

Table~\ref{tab1} is not meant to show that \pcp is clearly superior to all other
alternatives on all or most instances. While we aim at being very competitive in speed,
this work is not a competition paper; we find such quest 
quite absurd. The three
compared algorithms rely on different philosophies. The speed of \pcp
depends on the way \matlab implements certain basic operations (like Cholesky or QR 
factorizations,
back-substitution, matrix multiplication, etc). 
Most of these building blocks have a theoretical complexity of $O(n^3)$
but their running time in \matlab (version 2018) seem closer to $O(n^2)$.
This explains why the last column of Table~\ref{tab1} suggest
that \pcp is the most competitive method for $n\geq 2000$.
%(which is corroborated with the results from Table~\ref{tab2}).

\begin{table}
\hspace{-1.2em}
{\scriptsize
\begin{tabular}{cccc||c|c|c|c|c|c||c|c||c}
\multicolumn{4}{c||}{Instance} & 
\multicolumn{6}{c||}{\pcp} & 
\multicolumn{2}{c||}{\conb} &
\mosek \\
\drow{$n$} & \drow{k} &Eigs& Eigs  & Itera- & All          & Compute & Proj      & LP time & Send data & Trace   & Trace     &    \\
           &      & $A_i$'s&   C   & tions  &   time       &  $X$ \& $D$&time    & (cplex) &  to LP    & unknown & provided  &    \\
\hline
%n  & k  &  Eigs     & Eigs C & Iterations & Wall time & $X$ \& $D$ & Proj time & LP time & Send data to LP & Trace unknown & Trace known\\
800 &80  &[-20,~100]&[0,100] &      1108  & 410       &   179      &  44       &70       &   102           &  1051         & 94         & 320   \\
600 &40  &[-20,~100]&[0,100] &      155   & 17        &   4        &   6       &1        &   3             &  148          & 22         & 72 \\
400 &100 &[-20,~100]&[0,100] &      2075  & 572       &   94       &  13       & 384     &   71            &  490          & 42         & 60  \\
%200 &150 &[50,100] &[50,100]&      2767  & 1314      &   85       &  8        & 1149    &   56            & 6258          & 172        &    \\
%200 &200 &[20,21]  &[1,2]   & 3681       & 2866      &   154      &  28       &2556     &  107            & 14983         & 1734       &    \\
\hline
\multicolumn{12}{c}{Huge instances below have $\y\geq 0$, a random $\b$ and $\frac n5$ fixed null eigenvectors for all $A_i$'s and $C$}\\
\hline
200 &2000&[40, 100] &[10,40] &   31       &  11       &  5        & 0.2       & 0.2     &  5              &\multicolumn{2}{c||}{timeout}& 717     \\
200 &3000&[40, 100] &[10,40] &   70       &  49       &  27        & 0.4       & 0.7     &  18             &\multicolumn{2}{c||}{timeout}& 1346    \\
4000&20  &[20,25]   &[20,25] &  8         &76         &   17       & 44        &  0      &  11             &\multicolumn{2}{c||}{timeout}& timeout \\
5000&20  &[20,25]   &[20,25] &  7         &139        &   27       & 87        &  0      &  18             &\multicolumn{2}{c||}{timeout}& timeout \\

\end{tabular}
}
\caption{\small Seven runs of \pcp, \conb and \mosek on more varied instances. The last four instances
have $\y\geq \zeros$; such linear constraints on $\y$ simplify the
problem for \pcp, but this may be a non-trivial change for \conb (or other algorithms
that do not embed the SDP problem in a lightweight LP over $\y$).\label{tab2}}
\end{table}

Table~\ref{tab2} next page compares \projcutplanes with the \conb on the second benchmark 
set with instances of more varied sizes and
of a different nature (regarding the spectrum of the $A_i$'s or the non-negativity of $\y$).
Switching to $\y\geq 0$ may heavily reduce the number of \pcp iterations because most of
the elements of the optimal 
$\y$ may be zero.
In some cases, even if $k$ reaches a value of thousands, the associated LPs
remain very easy in practice because many of the $\y$ variables may remain zero at optimality when
$\y\geq \zeros$.
The last four rows of this table suggest \pcp is the best method for very large
SDP programs.\footnote{We provided the optimal trace to the \conb in the 
run from the last column.
Since the optimal trace is unknown in advance, we determined it from the \conb run from
the next-to-last column where we only used a bounded trace constraint 
(as in Remark~\ref{rem3}, p.~\pageref{rem3}). We are fully conscious that a better 
implementation of this optimal trace constraint may speed up the \conb.}

We presented up to here only the most relevant benchmarking information we could
present in a 12-pages paper. But the results reported in this work-in-progress article are not
a perfect measure of the final potential of the projection 
idea. 
This work represens the most initial version of the proposed method, submitted 
for the very first time to peer review. We 
must confess such software can not be perfect,
because it was not thoroughly tested. Perhaps other SDP 
algorithms
out there invested 
1000 times more coding and software testing resources.
However, it is quite safe and easy to check the correctness of a lower 
bound $\b^\top \overline\y$ reported by \pcp: it is enough to check that the minimum eigenvalue 
of $C-\AA^\top \overline\y$ is not-negative. It is very
difficult to have errors in the upper bounds either, because any $\d\in\R^n$ provides
a valid cut \eqref{sdp2d} and each upper bound is simply computed by the
LP solver that optimizes over all cuts \eqref{sdp2d} provided all along the iterations.

\section{Conclusion and prospects\label{sec5}}

We used \pcp ideas~\cite{mainpaper} to propose a fast method for optimizing (very)
large SDP programs.
Many ideas go beyond SDP optimization, because the
considered SDP program is incorporated in a more general (and yet very simple) LP.
For example, the \cutplanes logic for solving this LP enables one to easily insert some 
initial linear constraints \eqref{sdp2c} in the main SDP problem
    \eqref{sdp2a}-\eqref{sdp2d}.
If these linear constraints are prohibitively-many, they could even be generated on the fly by solving a second projection sub-problem in 
a purely LP context.
    We plan to implement this idea on a robust SDP problem in which
    the coefficients of the nominal constraints \eqref{sdp2c} can vary according to some robust rules and produce prohibitively-many 
     robust cuts -- a projection algorithm for
    this robust LP is already available~\cite[Sec.~2.1]{followup}.
It is quite easy to adapt \pcp to perform certain re-optimization tasks like the following:
after solving a \eqref{sdp2a}-\eqref{sdp2d} program, solve the same program again
after adding a new LP (or SDP) constraint.
    We are not aware of other methods that can adapt so easily to
    address such questions.

\bibliographystyle{splncs04}
\bibliography{all}
\appendix
\section{More insights into the design and the implementation of \pcp}

While the key element of this 
work is the projection sub-problem, the overall implementation 
depend on many other (down-to-earth) factors. The main paper presented only 
the most important guidelines for understanding \pcp, but it is 
not possible to discuss all nuts and bolts of the method. This 
appendix provides a number of remarks that completes the 
description of certain components of \pcp.

\begin{remark}\label{remConv}
In theory, the feasible area $\PP$ of \eqref{sdp2a}-\eqref{sdp2d} is not a polytope. But if we consider
in \eqref{sdp2d} 
only constraints $\d$ with a finite number of digits, this feasible area becomes a 
polytope.
As long as the amount of memory available on Earth is finite, the infinite
number of SDP cuts is actually finite when one solves \eqref{sdp2a}--\eqref{sdp2d} with an earthly 
computer. A \pcp iteration $\itt$ either returns a new cut never discovered before or
stops by proving $\optsol (\PP_{\itt})$ is optimal (with $t^*=1$). Considering a finite
number of potential SDP cuts, the algorithm will converge in a finite number of
iterations. We could go into more technical questions on convergence proofs, but 
such techniques are not directly relevant to the core of our algorithms. \qed
\end{remark}

\begin{remark}\label{remCaseD}
Projecting $X\to D$ is equivalent to projecting $\midx\to \midmat$ using
\eqref{eqinner2} as discussed at point D) of the projection algorithm 
(p.~\pageref{eqinner2}).
It may be faster to use the smaller matrices $\midx$ and $\midm$ of order $c+m<n$.
On the other hand, not working with the original matrices may lead to more numerical 
problems. We here limit the presentation to the case in
which we apply repeated separation on the original matrices $X$ and $D$: we have
to use repeated separation
to determine the SDP
status of $X+tD$ for various values of $t$.

We consider a user-provided list of separation points $t_1,~t_2,~t_3, \dots$ to be tried
so that $0<t_1<t_2<t_3\dots$; for each such $t_i$, we solve the separation sub-problem by 
determining the minimum eigenvalue of  $X+t_i D$. If this value is negative, $X+t_i D$ does
not belong to the SDP cone. We now split case D) in two sub-cases:
\begin{description}
\item[D.1)] If $X+t_1D\nsucceq \zeros$ we return $t^*=0$. It is important to have 
a $t_1$ value (very) close to 0. We basically consider that there is no space inside the
SDP cone to perform any positive step towards $D$ only because there is no space to
perform a step of $t_1$. By using a $t_1$ close to $10^{-6}$ we avoid many numerical 
problems when solving case D.1) this way.
\item[D.2)] If $X+t_1D\succeq \zeros$, we return $t^*=t_1+t^*_2$, where $t^*_2$ is the
step length returned by projecting $(X+t_1D)\to D$. But the advantage of this new  
projection is that $D$ will belong to the image of $X+t_1D$, except in very pathological 
cases.
This way, Property~\ref{prop1} is very likely to hold and we can solve the projection
using cases A) or B).\qed
\end{description}

\end{remark}

\begin{remark}\label{rem3}
To provide a constant trace constraint for the \conb in the dual of
\eqref{sdp2a}-\eqref{sdp2d},
we insert into the primal \eqref{sdp2a}-\eqref{sdp2d} an additional $A_{k+1}=I_n$ 
alongside a $b_{k+1}$ equal
to the trace value. 
Since the trace is unknown in our experiments, we can only provide it by solving
the instance beforehand.
\pcp could have also exploited 
such information to produce more interior feasible solutions.
However, in all \conb experiments with an unknown trace, we indicated to \conb that this trace is 
$\leq 1000$. For this, we added a row and column
of zeros to all matrices $A_1,A_2,\dots A_{k+1}$ and $C$, putting a 1 only at position $(n+1,n+1)$
of $A_{k+1}$ so that $A_{k+1}$ becomes $I_{n+1}$, and $b_{k+1}=1000$. 
We are fully conscious that a better 
implementation of this optimal trace constraint may speed up the \conb.
\qed
\end{remark}

\begin{remark}\label{remNutsAndBolts} The most important customizations of
\pcp that were not fully described in the main body of the paper (due to space limitation)
are the
following.
\begin{itemize}
\item In the very beginning there is no default constraint that \pcp may  
use to construct a very first outer approximation of \eqref{sdp2a}-\eqref{sdp2d} or
a very first outer solution.
We inserted an artificial initial box to have such a first outer approximation. This box
only limits each variable $y_i$ to the interval [-100000,100000], which is more than 
enough for our instances.
In the very beginning, while the current $\y$ still touches the box,
 we use standard \cutplanes (this never took an important amount of time compared
 to overall running time).

\item Recall that in cases A) and B) we computed the minimum eigenvalue of $D'$ 
in \eqref{eqmaxt3}, or respectively, \eqref{eqmaxtb}. We described how
that minimum eigenvalue
produces a first-hit cut. We noticed that in practice it may be
useful to go to the second minimum eigenvalue and use it to compute a second-hit cut using
exactly the same calculations as for the first-hit cut. We certainly 
do this only if this second minimum eigenvalue is still negative.

\item We normalize certain cuts we eventually send to the LP
solver ({\tt cplex}). 
For each $i\in[1..k]$ the coefficient $i$ of decision variable $y_i$ comes from
the term $\v^\top A_i \v$, where $\v$ is the first-hit vector returned
by the projection algorithm.
When the maximum resulting coefficient in absolute value is greater than 
100000, we divide all coefficients of the cut by that maximum
coefficient.
\qed 
\end{itemize}
\end{remark}

\begin{remark}\label{remComputer}
The code was implemented in \matlab (version {\tt r2018b}) on a mainstream
laptop clocked at 1.90GHz with an Intel i7-8665U processor with 4 cores. The number
of threads can go up to 8 using a hyper-threading technology. We used the default
\matlab configuration that allows up to 4 threads (a {\tt maxNumCompThreads} value
of 4). We chose \matlab because  preliminary experiments suggest it provides the fastest 
matrix eigenvalue
routines for $n\geq 1000$.
We used a Linux Mint operation system; the Linux kernel version is 4.15.0. The LP solver is 
{\tt cplex} version 12.10. \qed
\end{remark}

\end{document}

%% file: figures/fig1.pdf_tex
%% Creator: Inkscape inkscape 0.91, www.inkscape.org
%% Accompanies image file 'fig1.pdf' (pdf, eps, ps)
\begingroup%
  \makeatletter%
  \providecommand\color[2][]{%
    \errmessage{(Inkscape) Color is used for the text in Inkscape, but the package 'color.sty' is not loaded}%
    \renewcommand\color[2][]{}%
  }%
  \providecommand\transparent[1]{%
    \errmessage{(Inkscape) Transparency is used (non-zero) for the text in Inkscape, but the package 'transparent.sty' is not loaded}%
    \renewcommand\transparent[1]{}%
  }%
  \providecommand\rotatebox[2]{#2}%
  \ifx\svgwidth\undefined%
    \setlength{\unitlength}{512.46154785bp}%
    \ifx\svgscale\undefined%
      \relax%
    \else%
      \setlength{\unitlength}{\unitlength * \real{\svgscale}}%
    \fi%
  \else%
    \setlength{\unitlength}{\svgwidth}%
  \fi%
  \global\let\svgwidth\undefined%
  \global\let\svgscale\undefined%
  \makeatother%
  \begin{picture}(1,0.39609364)%
    \put(0,0){\includegraphics[width=\unitlength,page=1]{fig1.pdf}}%
    \put(0.02322072,0.01828722){\color[rgb]{0,0,0}\makebox(0,0)[lb]{\smash{$\x_1=[0~0]^\top$}}}%
    \put(0,0){\includegraphics[width=\unitlength,page=2]{fig1.pdf}}%
    \put(0.14105755,0.15059786){\color[rgb]{0,0,0}\rotatebox{45.67525057}{\makebox(0,0)[lb]{\smash{\scalebox{0.95}{towards $\optsol(\PP_0)$}}}}}%
    \put(0,0){\includegraphics[width=\unitlength,page=3]{fig1.pdf}}%
    \put(0.02311246,0.03910929){\color[rgb]{0,0,0}\rotatebox{45.24519654}{\makebox(0,0)[lb]{\smash{iteration 1}}}}%
    \put(0.16884111,0.11410123){\color[rgb]{0,0,1}\rotatebox{-8.12489944}{\makebox(0,0)[lb]{\smash{iteration 2}}}}%
    \put(0.28529123,1.49200605){\color[rgb]{0,0,0}\makebox(0,0)[lt]{\begin{minipage}{0.28336704\unitlength}\raggedright \end{minipage}}}%
    \put(0.34308509,0.11166409){\color[rgb]{1,0,0}\rotatebox{20.12201124}{\makebox(0,0)[lb]{\smash{iteration 3}}}}%
    \put(0.02524884,0.30735549){\color[rgb]{0,0,0}\makebox(0,0)[lt]{\begin{minipage}{0.17285388\unitlength}\raggedright \scalebox{1.7}{$\PP$}\end{minipage}}}%
    \put(0.88948872,0.06198597){\color[rgb]{0,0,0}\makebox(0,0)[lt]{\begin{minipage}{0.06159942\unitlength}\raggedright \end{minipage}}}%
    \put(0.88948872,0.05098608){\color[rgb]{0,0,0}\makebox(0,0)[lt]{\begin{minipage}{0.05829953\unitlength}\raggedright \end{minipage}}}%
    \put(0.91184012,0.04075582){\color[rgb]{0,0,0}\makebox(0,0)[lt]{\begin{minipage}{0.26591801\unitlength}\raggedright $\optsol(\PP_1)$\end{minipage}}}%
    \put(0.10443499,0.14685164){\color[rgb]{0,0,1}\makebox(0,0)[lt]{\begin{minipage}{0.07793552\unitlength}\raggedright $\x_2$\end{minipage}}}%
    \put(0.275363,0.10319686){\color[rgb]{0,0,0}\makebox(0,0)[lt]{\begin{minipage}{0.09715211\unitlength}\raggedright \end{minipage}}}%
    \put(0.2888607,0.10650191){\color[rgb]{1,0,0}\makebox(0,0)[lt]{\begin{minipage}{0.10531349\unitlength}\raggedright $\x_3$\end{minipage}}}%
    \put(0,0){\includegraphics[width=\unitlength,page=4]{fig1.pdf}}%
    \put(0.28498332,1.59255888){\color[rgb]{0,0,0}\makebox(0,0)[lt]{\begin{minipage}{0.07124649\unitlength}\raggedright \end{minipage}}}%
    \put(0.25742714,0.2993368){\color[rgb]{0,0,0}\makebox(0,0)[lt]{\begin{minipage}{0.0831844\unitlength}\raggedright \end{minipage}}}%
    \put(0.27215902,0.28523992){\color[rgb]{0,0,0}\makebox(0,0)[lt]{\begin{minipage}{0.08991535\unitlength}\raggedright $\x_1+t^*_1\d_1$\end{minipage}}}%
    \put(0.38325764,0.08382321){\color[rgb]{0,0,1}\makebox(0,0)[lt]{\begin{minipage}{0.08991535\unitlength}\raggedright $\x_2+t^*_2\d_2$\end{minipage}}}%
    \put(0,0){\includegraphics[width=\unitlength,page=5]{fig1.pdf}}%
    \put(0.34650106,0.17937523){\color[rgb]{1,0,0}\makebox(0,0)[lt]{\begin{minipage}{0.08991535\unitlength}\raggedright $\x_3+t^*_3\d_3$\end{minipage}}}%
    \put(0.49951862,0.19924126){\color[rgb]{0,0,1}\makebox(0,0)[lt]{\begin{minipage}{0.26591801\unitlength}\raggedright $\optsol(\PP_2)$\end{minipage}}}%
    \put(0,0){\includegraphics[width=\unitlength,page=6]{fig1.pdf}}%
  \end{picture}%
\endgroup%